\documentclass{amsart}

\newtheorem{theorem}{Theorem}[section]

\newtheorem{proposition}[theorem]{Proposition}

\theoremstyle{definition}

\newtheorem*{acknowledgement}{Acknowledgement}
\newtheorem{question}[theorem]{Question}

\theoremstyle{remark}

\newcommand\mynote[1]{\marginpar{\ \\ \small \tt #1}}
\newcommand\bel[1]{{\mynote{#1}}\begin{equation}\label{#1}}
\newcommand\mylabel[1]{\label{#1}}

\newcommand{\ZZ}{\mathbb{Z}}

\newcommand{\GG}{\mathbb{G}}

\newcommand  {\shC}     {\mathcal{C}}

\newcommand  {\shF}     {\mathcal{F}}

\newcommand  {\shM}     {\mathcal{M}}

\newcommand  {\shL}     {\mathcal{L}}

\newcommand  {\Br}      {\operatorname{Br}}

\newcommand  {\Div}     {\operatorname{Div}}

\newcommand  {\et}      {{\text{\rm \'{e}t}}}
\newcommand  {\Et}      {{\text{\'{E}t}}}

\newcommand  {\Hom}     {\operatorname{Hom}}

\newcommand  {\id}      {\operatorname{id}}

\newcommand  {\lra}     {\longrightarrow}

\renewcommand{\O}       {\mathcal{O}}

\newcommand  {\Pic}     {\operatorname{Pic}}

\newcommand  {\quadand} {\quad\text{and}\quad}

\newcommand  {\ra}      {\rightarrow}

\newcommand  {\Sch}     {\operatorname{Sch}}

\newcommand  {\Spec}    {\operatorname{Spec}}

\newcommand  {\Supp}    {\operatorname{Supp}}

\newcommand {\zar}      {{\operatorname{zar}}}

\def\mydate{\number\day\space\ifcase\month \or January\or February\or March\or 
April\or May\or June\or July\or
August\or September\or October\or November\or December\fi \space\number\year}
\usepackage{amscd,amssymb}

\begin{document}

\title[Hilbert's Theorem 90 and algebraic spaces]
      {Hilbert's Theorem 90 and algebraic spaces}

\author[Stefan Schroer]{Stefan Schr\"oer}

\address{Fakult\"at f\"ur Mathematik, Ruhr-Universit\"at, 
         44780 Bochum, Germany}

\email{s.schroeer@ruhr-uni-bochum.de}

\subjclass{14A20, 14C22, 14F20}

\dedicatory{Revised version, 22 Oktober 2001}

\begin{abstract}
In modern form, 
Hilbert's Theorem 90 tells us that $R^1\epsilon_*(\GG_m)=0$, where
$\epsilon:X_\et\ra X_\zar$ is the canonical map between 
the \'etale site and the Zariski site of a scheme $X$.
I construct examples showing that the corresponding statement
for algebraic spaces does not hold. 
\end{abstract}

\maketitle

%===========================================================
\section*{Introduction}

Originally, Hilbert's Theorem 90  is the
following number theoretical result \cite{Hilbert 1897}:
Given a cyclic Galois extension $K\subset L$ of number fields,
each $y\in L^\times$ of norm $N(y)=1$ is of the form $y=x/x^\sigma$ for some
$x\in K^\times$ and a given generator $\sigma\in G$ of the Galois group.
More generally, Speiser \cite{Speiser 1919} proved that
$H^1(G,L^\times)=1$ for arbitrary Galois extensions
(compare the discussion in \cite{Lorenz 1998}).

The latter statement has a geometric interpretation:
Each line bundle on the \'etale site of $\Spec(k)$ is trivial.
In this form, it admits a far-reaching  generalization:
If  $\epsilon:X_\et\ra X_\zar$ is the canonical map from the
\'etale site to the Zariski site of a scheme $X$, then 
$R^1\epsilon_*(\GG_m)=0$ (see \cite{Milne 1980}, page 124).
The result entails, among other things, that the map of Picard groups
$\Pic(X_\zar)\ra \Pic(X_\et)$ is bijective, 
and that the map of Brauer groups $\Br(X_\zar)\ra\Br(X_\et)$ is injective.

It is natural to ask whether a similar statement holds for 
algebraic spaces instead of schemes. Recall that an \emph{algebraic space}
is the quotient $X=U/R$ of a scheme $X$
by an \'etale equivalence relation $R\rightrightarrows X$.
Here the quotient takes place in the topos $(\Sch)_\et^\sim$, that
is, as a sheaf on the \'etale site.

Unfortunately, such a generalization does not hold.
The goal of this paper is to construct counterexamples, that is,
algebraic spaces $X$ and  invertible $\O_X$-modules
$\shL$ such that the open subspaces $V\subset X$ trivializing $\shL$
do not cover $X$.
The first example is a nonseparated smooth 1-dimensional
\emph{bug-eyed cover} in Koll\'ar's sense \cite{Kollar 1992}.
The second example is a nonnormal proper algebraic space obtained by 
identifying  points on suitable nonprojective smooth proper schemes.

\begin{acknowledgement}
The author wishes to thank the Department of Mathematics of the
Massachusetts Institute of Technology for its hospitality,
and the Deutsche For\-schungs\-gemein\-schaft for financial support.
The author also thanks the referee for his interesting remarks.
\end{acknowledgement}

%===========================================================
\section{Line bundles on algebraic spaces}

In this section we recall some basic facts on algebraic spaces and
their line bundles.
Let $(\Sch)_\et$ be the site  of schemes endowed with the
Grothendieck topology generated by the \'etale surjective morphisms,
and $(\Sch)_\et^\sim$ be the corresponding topos of sheaves.
By definition, a sheaf 
$X\in(\Sch)_\et^\sim$ is an \emph{algebraic space} if
$X=U/R$ for some scheme $U$ and some \'etale equivalence relation
$R\rightrightarrows U$ such that the induced morphism  $R\ra U\times U$ 
is quasicompact
\cite{Knutson 1971}.

Given an algebraic space $X$, let $\Et(X)$ be the category of
algebraic $X$-spaces whose structure map $Y\ra X$ is \'etale.
The \'etale surjections $Y_1\ra Y_2$ define a topology on $\Et(X)$, 
and we write
$X_\et$ for the corresponding site.
Let me give a down-to-earth description of sheaves $\shF$ on this site.
For each  scheme $U$ endowed with an \'etale map $U\ra X$,
we obtain via restriction a sheaf $\shF_U$ on the \'etale site
of  \'etale $U$-schemes. If $f:U\ra V$ is an $X$-morphisms,
we have a map $\theta_f:\shF_V\ra f_*\shF_U$.
Such systems $(\shF_U,\theta_f)$ are not arbitrary. 
Consider  the following two conditions:
(1) If $f:U\ra V$ and $g:V\ra W$ are $X$-maps, then the diagram
$$
\begin{CD}
\shF_W @>\theta_{gf}>> (gf)_*\shF_U\\
@V\theta_gVV @VV\simeq V\\
g_*(\shF_V) @>>g_*(\theta_f)> g_*(f_*\shF_U)
\end{CD}
$$
is commutative. (2) If $f:U\ra V$ is \'etale, then the map
$\theta_f^\sharp:f^{-1}\shF_V\ra\shF_U$
is bijective. Here the mapping 
$\theta_f^\sharp$ corresponds to $\theta_f$ with respect to  the canonical 
adjunction
$\Hom(f^{-1}\shF_V,\shF_U)\simeq\Hom(\shF_V,f_*\shF_U)$.

\begin{proposition}
\mylabel{equivalence}
The assignment $\shF\mapsto(\shF_U,\theta_f)$ 
yields an equivalence between the category of sheaves on $X_\et$ and
the category of systems $(\shF_U,\theta_f)$ satisfying conditions 
{\rm (1)} and {\rm (2)}.
\end{proposition}

\proof
Let $\shC$ be the site of \'etale $X$-schemes with the 
induced \'etale topology. By the Comparison Lemma
(\cite{SGA 4a}, Expos\'e III, Th\'eor\`em 4.1),
the inclusion $\shC\subset X_\et$ induces an equivalence on 
the corresponding categories of sheaves.
Now suppose  $\shF$ is a sheaf on $\shC$.
Then the system $(\shF_U,\theta_f)$ satisfies condition (1) because 
$\shF$ is a 
presheaf.
If $f:U\ra V$ is \'etale, then $\theta_f^\sharp$ is bijective because 
$\shF$ is a sheaf in the \'etale topology, and condition (2) holds
as well.

Conversely, given such a system, we define
$\Gamma(U,\shF)=\Gamma(U,\shF_U)$. 
Indeed, this is a presheaf by condition (1), and a sheaf
by condition (2). One easily checks that
the functors $\shF\mapsto(\shF_U,\theta_f)$ and 
$(\shF_U,\theta_f)\mapsto\shF$ are inverse equivalences of categories.
\qed

\medskip
For example, the sheaves $\O_U$, together with the  maps
$\theta_f:\O_V\ra f_*(\O_U)$, 
correspond to the structure sheaf $\O_X$ of an algebraic space $X$.
Similar, we have the sheaf of units $\O_X^\times$.
The cohomology group $\Pic(X_\et)=H^1(X_\et,\O_X^\times)$ 
is the group of isomorphism classes of 
invertible $\O_X$-modules.

Besides the \'etale topology, the 
category $\Et(X)$  carries the coarser \emph{Zariski topology}
as well. Here the 
covering families are the surjections of the form $\coprod X_i\ra X$,
where the $X_i\subset X$ are open subspaces, and we demand that
$X_i\times_X X'\ra X'$ remains an open embedding for any base
change
$X'\ra X$. Write $X_\zar$ for the corresponding site.
The sheaves on $X_\zar$ admit a similar description in terms of families
$(\shF_U,\theta_f)$ satisfying condition (1), and condition (2'), where
we demand that $\theta_f^\sharp:f^{-1}\shF_V\ra \shF_U$ 
is bijective whenever $f:U\ra V$ is of the form $U=\coprod V_i$
with open subschemes $V_i\subset V$.
In particular, we have a structure sheaf $\O_{X_\zar}$ and a unit sheaf 
$\O_{X_\zar}^\times$.
Let $\Pic(X_\zar)=H^1(X_\zar,\O_{X_\zar}^\times)$ 
be the corresponding group of line bundles.

The identity functor on $\Et(X)$ is a continuous functor 
$\epsilon:X_\et\ra X_\zar$
of sites, and we have $\epsilon_*(\O_{X_\et})=\O_{X_\zar}$ by
descent theory. So for each invertible $\O_{X_\zar}$-module $\shL$, 
the canonical map $\shL\ra\epsilon_*\epsilon^*\shL$ 
is bijective, and we obtain an
injection $\Pic(X_\zar)\subset \Pic(X_\et)$.

\begin{proposition}
\mylabel{line bundle}
Let $\shL$ be an invertible $\O_X$-module.
Its isomorphism class lies  in the subgroup
$\Pic(X_\zar)\subset\Pic(X_\et)$ 
if and only if there is a covering with open subspaces
$Y_i\subset Y$ with  $\shL_{Y_i}\simeq \O_{Y_i}$.
\end{proposition}

\proof
The spectral sequence for the composition 
$\Gamma(X_\et,\O_{X_\et}^\times)=\Gamma(X_\zar,\epsilon_*\O_{X_\et}^\times)$
yields an exact sequence
$$
0\lra \Pic(X_\zar)\lra\Pic(X_\et)\lra 
H^0(X_\zar,R^1\epsilon_*\O_{X_\et}^\times).
$$
The condition precisely means that the image of the invertible sheaf $\shL$
under the canonical 
map $\Pic(X_\et)\ra H^0(X_\zar,R^1\epsilon_*\O_{X_\et}^\times)$ vanishes. 
The statement now follows from the exact sequence.
\qed

%===========================================================
\section{Bug-eyed covers}

In this section, we use Koll\'ar's bug-eyed covers to construct
a smooth 1-dimen\-sional \emph{nonseparated}
algebraic space $X$ and an invertible sheaf 
$\shL$ such that the open 
subspaces $W\subset X$ trivializing $\shL$ do not form a covering.

Fix a ground field $k$ of characteristic $\neq 2$.
Set $A=k[[T]]$ and $A'=k[[T^2]]$, and let 
$Y=\Spec(A)$ and $Y'=\Spec(A')$ be the corresponding affine schemes.
The inclusion $A'\subset A$ defines a flat double covering $p:Y\ra Y'$.
The open subset $U\subset Y$ given by the generic point is the locus
where $f$ is \'etale. The generator $\sigma\in G$ of the group $G=\ZZ/2\ZZ$
acts on $A$ via $T^\sigma=-T$, which defines a free $G$-action on $U$.
Consider the \'etale equivalence relation 
$$
R=\Delta_Y\amalg U\lra Y\times Y,
$$
where the embedding of $U$ is given by
$U\stackrel{\id\times \sigma}{\lra}U\times U\subset Y\times Y$.
Let $X=Y/R$ be the corresponding quotient sheaf in $(\Sch/k)_\et^\sim$.
By definition, $X$ is a smooth algebraic space. It is nonseparated
because the injection $R\ra Y\times Y$ is not closed.

The map $p:Y\ra Y'$
factors over $X$, and the induced projection $X\ra Y'$ induces
a bijection of points.
The algebraic space $X$ is a \emph{bug-eyed cover}
in Koll\'ar's sense \cite{Kollar 1992}.
It is not a scheme. Otherwise, the morphism $X\ra Y'$ 
would be an isomorphism
by Zariski's Main Theorem, and $Y\ra X$ would be both \'etale and ramified.

\begin{proposition}
\mylabel{picard}
We have $\Pic(X_\et)=\ZZ/2\ZZ$.
\end{proposition}

\proof
The scheme $Y$ is local, hence every invertible $\O_X$-module
$\shL$ has $\shL_Y\simeq\O_Y$. Thus, $\Pic(X_\et)$ 
is the cohomology of the complex
$$
\Gamma(Y,\O_X^\times)\stackrel{d_0}{\lra}
\Gamma(Y^2,\O_X^\times)\stackrel{d_1}{\lra}
\Gamma(Y^3,\O_X^\times).
$$
Here $Y^n$ are the $n$-fold fiber products over $X$. 
If $p_i:Y^{n+1}\ra Y^n$ 
denotes the projection omitting the $i$-th factor, the
differentials are $d_0(s)=p_0^*(s)/p_1^*(s)$ and 
$d_1(s)=p_0^*(s)p_2^*(s)/p_1^*(s)$.

Clearly, we have $Y^n=U^n\cup\Delta_Y$, where $U^n\cap\Delta_Y=\Delta_U$.
Since the $G$-action is free on the open subset $U\subset Y$, we have a
bijection
$$
U\times G^n\lra U^{n+1},\quad 
(u,g_1,\ldots,g_n)\longmapsto (u,ug_1,\ldots,ug_1g_2\ldots g_n).
$$
In turn, we may identify the $n$-cochains $\Gamma(Y^{n+1},\O_X^\times)$ 
with the
the group of functions $c:G^{n}\ra P^\times$ satisfying 
$c(0,\ldots,0)\in A^\times$.
Here $P=k[[T]][T^{-1}]$ is the fraction field of $A=k[[T]]$.
The differentials take the form
$$
d_0(c)(g)=c(0)/c(0)^g \quadand
d_1(c)(g,h)=c(h)^g c(g)/c(gh),
$$
conforming with the usual definition of group cohomology
(\cite{Brown 1982}, page 59).
We have $d_0(c)(0)=1$, and $d_0(c)(\sigma)$ is a power series of the form
$\lambda_0+\lambda_1T+\lambda_2T^2+\ldots$ with $\lambda_0=1$. 
One easily checks that a 1-cochain $c:G\ra P^\times$ is a 1-cocycle
if and only if $c(0)=1$, and $p=c(\sigma)$ satisfies $p\cdot p^\sigma=1$.
Clearly, the 1-cocycle $c:G\ra P^\times$ with
$c(0)=1$ and $c(\sigma)=-1$ is not a coboundary, so $\Pic(X_\et)$ is nonzero.
On the other hand, by Hilbert's Theorem 90, each $p\in P^\times$ 
with $p\cdot p^\sigma=1$
is of the form $p=r/r^\sigma$ for some $r\in P^\times$.
Writing $r=T^n s$ with $s\in A^\times$, we have $p=(-1)^n s/s^g$, 
and infer $\Pic(X_\et)=\ZZ/2\ZZ$. 
\qed

\medskip
The smooth 1-dimensional nonseparated algebraic space $X$ is 
our first counterexample to Hilbert's Theorem 90 for algebraic spaces:

\begin{theorem}
\mylabel{nonzero1}
The canonical inclusion $\Pic(X_\zar)\subset\Pic(X_\et)$ is not surjective.
\end{theorem}

\proof
The scheme $Y$ is local, so the space of points for $X$ has a unique
closed point. Consequently, any Zariski covering of $X$ contains
a copy of $X$. So any line bundle on
$X_\zar$ is trivial, that is, $\Pic(X_\zar)=0$. 
On the other hand, $\Pic(X_\et)\neq 0$ by Proposition \ref{picard}.
\qed

%===========================================================
\section{Nonnormal proper algebraic spaces}

Fix an algebraically closed ground field $k$.
In this section, we shall construct a  \emph{proper}  algebraic space
$X$ and an invertible sheaf $\shL$ such that the open 
subspaces $W\subset X$ trivializing $\shL$ do not form a covering.

The starting point is a proper smooth $k$-scheme $Y$
containing two irreducible closed  curves $C_1,C_2\subset Y$ 
such that $C_1+C_2$ is 
numerically trivial. 
This implies that the generic points
$\eta_i\in C_i$ do not admit any common affine
neighborhood in $Y$. Examples of such schemes appear
in \cite{Shafarevich 1994}, page 75. Obviously, they are nonprojective.
Even worse, they do not admit embeddings into
toric varieties (\cite{Wlodarczyk 1993}, Theorem A).
Recall that the support $\Supp(D)\subset Y$ of a Cartier divisor $D\in \Div(Y)$
is the union of its positive and negative part.
We have the following useful property:

\begin{proposition}
\mylabel{cartier}
Each  $D\in\Div(Y)$ with $D\cdot C_1>0$ and
$C_1\not\subset\Supp(D)$ has 
$C_2\subset\Supp(D)$.
\end{proposition}

\proof
Decompose $D=\sum n_iD_i$ into prime divisors with $n_i\neq 0$. 
Since $C_1\not\subset D_i$, the intersection number
$D_i\cdot D_1$ is the length of the scheme $D_i\cap C_1$, 
hence nonnegative. So there is at least one prime divisor
with $D_i\cdot C_1>0$.
It follows $D_i\cdot C_2<0$, hence $C_2\subset D_i$.
In other words, $C_2\subset\Supp(D)$.
\qed

\medskip
Now fix two closed points $y_1\in C_1$ and $y_2\in C_2$.
Let $Y'\subset Y$ be the reduced closed subscheme corresponding to 
$\left\{y_1,y_2\right\}$, and
define an \'etale sheaf $X\in(\Sch/k)_\et^\sim$ by the cocartesian square
$$
\begin{CD}
Y' @>>> Y\\
@VVV @VVp V \\
\Spec(k) @>>> X
\end{CD}
$$
Note that $(\Sch/k)_\et^\sim$, being a topos, admits all colimits
(\cite{SGA 4a}, Expos\'e II, Theorem 4.1).
Intuitively, $X$ is obtained from $Y$ by identifying the points
$y_1,y_2\in Y$.
The sheaf $X$ is not a scheme. Otherwise, an affine open neighborhood
for the point
$p(y_1)=p(y_2)\in X$ would give a common affine open neighborhood 
for the pair
$y_1,y_2\in Y$.

\begin{proposition}
\mylabel{algebraic space}
The \'etale sheaf $X$ is a proper algebraic space.
\end{proposition}

\proof
That $X$ is  an algebraic space
follows immediately from \cite{Artin 1970}, Theorem 6.1.
Let me give a more direct argument as follows.
Fix two copies $v_1',v_2'\in V'$ and $v_1'',v_2''\in V''$ of $y_1,y_2\in Y$,
and set $V=V'\amalg V''$.
Identifying $v_1'\in V$ with $v_2''\in V$
and $v_2'\in V$ with $v_1''\in V$, 
we obtain a scheme $U$. The group $G=\ZZ/2\ZZ$ acts freely on
$U$ by interchanging $V'$ and $V''$. Clearly, $X=U/G$ is the quotient of
this action in the topos of \'etale sheaves.
So $R=U\times_X U$ is nothing but $U\times G$, which is a scheme.
Consequently, $X=U/R$ is an algebraic space.

The algebraic space $X$ is separated
because  the embedding $Y\times G\ra Y\times Y$,
$(y,g)\mapsto(y,yg)$ is closed.
As $Y\ra\Spec(k)$ is universally closed and $p:Y\ra X$ 
is surjective, $X\ra \Spec(k)$ is universally closed
as well. Therefore, $X$ is proper.
\qed

\begin{proposition}
\mylabel{surjective}
There is an exact sequence $1\ra k^\times\ra\Pic(X_\et)\ra\Pic(Y)\ra 0$.
\end{proposition}

\proof
Let $p:Y\ra X$ be the canonical projection.
Then the sequence
$$
1\lra \O_X^\times\lra p_*(\O_Y^\times)\oplus k^\times\lra
p_*(\O_{Y'}^\times)\lra 1
$$
is exact. Indeed, one easily checks this,
as in \cite{Hartshorne 1994}, Lemma 5.1, after base change
with an affine \'etale cover $U\ra X$.
In turn, we obtain an exact sequence
$$
\Gamma(\O_Y^\times)\oplus k^\times\lra\Gamma(\O_{Y'}^\times)
\lra\Pic(X_\et)\lra \Pic(Y)\oplus\Pic(k)\lra\Pic(Y').
$$
Being semilocal, the schemes $\Spec(k)$ and $Y'$ have no Picard groups.
The cokernel for the map on the left is isomorphic to $k^\times$,
and the result follows.
\qed

\medskip
The proper algebraic space $X$ is another counterexample
to Hilbert's Theorem 90 for algebraic spaces:

\begin{theorem}
\mylabel{nonzero2}
The canonical inclusion $\Pic(X_\zar)\subset\Pic(X_\et)$ is not surjective.
\end{theorem}

\proof
Choose an invertible $\O_Y$-module 
$\shM$ with $\shM\cdot C_1>0$.
For example, $\shM$ could by the invertible sheaf corresponding
to the  reduced complement
of any affine open neighborhood for $y_1\in Y$.

Let $p:Y\ra X$ be the canonical map.
According to Proposition \ref{surjective}, there is an invertible
$\O_X$-module $\shL$ with $\shM=p^*(\shL)$. Suppose there is 
an open subset $W\subset X$ containing the point $p(y_1)=p(y_2)$ 
and  trivializing $\shL$.
Then $\shM$ is trivial on the open subscheme $p^{-1}(W)\subset Y$.
By \cite{Schroeer 2000}, Theorem 3.3, there is a Cartier divisor
$D\in\Div(X)$ representing $\shM$ 
with support disjoint from $y_1,y_2\in Y$.
In particular, $C_1$ and $C_2$ are not contained in $\Supp(D)$,
contradicting  Proposition \ref{cartier}.
\qed

\begin{question}
Does $\Pic(X_\zar)=\Pic(X_\et)$ at least hold for smooth proper
algebraic spaces? What about the case that $X$ is normal and proper?
\end{question}

%===========================================================

\end{document}